\input amstex
\input amsppt.sty   
\def\nmb#1#2{#2}         
\def\ign#1{}             

\mag1200
\hsize 160truemm
\vsize 240truemm
\redefine\o{\circ}

\define\la{\lambda}

\define\ph{\varphi}

\define\om{\omega}

\define\Ph{\Phi}

\define\Om{\Omega}
\redefine\i{^{-1}}
\define\row#1#2#3{#1_{#2},\ldots,#1_{#3}}
\define\x{\times}
\define\Lip{\operatorname{Lip}}

\define\ev{\operatorname{ev}}
\def\today{\ifcase\month\or
 January\or February\or March\or April\or May\or June\or
 July\or August\or September\or October\or November\or December\fi
 \space\number\day, \number\year}
\topmatter
\title Aspects of the Theory of Infinite Dimensional Manifolds  \endtitle
\author Andreas Kriegl \\ 
Peter W. Michor  \endauthor
\affil
Institut f\"ur Mathematik, Universit\"at Wien,\\
Strudlhofgasse 4, A-1090 Wien, Austria.
\endaffil
\address{Institut f\"ur Mathematik, Universit\"at Wien,
Strudlhofgasse 4, A-1090 Wien, Austria.}\endaddress
\date{September 2, 1989}\enddate
\thanks{This paper is in final form and no version of it will appear 
elsewhere}\endthanks
\subjclass{26E15, 26E05, 46G20, 58B10, 58B12, 58D05, 58D10, 
58D15}\endsubjclass
\keywords{Convenient vector spaces, manifolds of mappings, 
diffeomorphism groups, infinite dimensional Lie groups and Grassmann 
manifolds}\endkeywords
\abstract{The convenient setting for smooth 
mappings, holomorphic mappings, and real analytic mappings in 
infinite dimension is sketched. Infinite dimensional manifolds are 
discussed with special emphasis on smooth partitions of unity and 
tangent vectors as derivations. Manifolds of mappings and 
diffeomorphisms are treated. Finally the differential structure on the 
inductive limits of the groups $GL(n)$, $SO(n)$ and some of their 
homogeneus spaces is treated.}\endabstract
\endtopmatter
\document

\heading Introduction \endheading

The theory of infinite dimensional manifolds has already a long
history. In the fifties and sixties smooth manifolds modeled on
Banach spaces were investigated a lot. Here the starting point
were the investigations of Marston Morse on the index of geodesics
in Riemannian manifolds. He used Hilbert manifolds of curves in
a Riemannian manifold. 

Later on Fr\'echet
manifolds were investigated from the point of view of topology:
it was shown that under certain weak conditions they could be
embedded as open subsets in the model space.

Then starting with a seminal short paper of J\. Eells began the
investigation of manifolds of mappings. 

But all this became important for the mainstream of mathematics
when loop groups and their Lie algebras - the Kac-Moody Lie
algebras were used in Physics.

In this review paper we will try to find our own way through the field of
infinite dimensional manifolds, and we will concentrate on the
smooth manifolds, and on those which are not modeled on Banach
spaces or Hilbert spaces - although the latter are very
important as a technical mean to prove very important theorems
like those leading to exotic $\Bbb R^4$'s. 

The material presented in the later sections is from 
\cite{Kriegl-Michor, Foundations of Global Analysis}.

\heading Table of contents \endheading

\noindent 1. Calculus of smooth mappings \par 
\noindent 2. Calculus of holomorphic mappings \par
\noindent 3. Calculus of real analytic mappings \par
\noindent 4. Infinite dimensional manifolds \par
\noindent 5. Manifolds of mappings \par
\noindent 6. Manifolds for algebraic topology \par

\heading\nmb0{1}. Calculus of smooth mappings \endheading

\subheading{\nmb.{1.1}} The traditional differential calculus works 
well for finite dimensional vector spaces and for Banach spaces. For 
more general locally convex spaces a whole flock of different 
theories were developed, each of them rather complicated and none 
really convincing. The main difficulty is that the composition of 
linear mappings stops to be jointly continuous at the level of Banach 
spaces, for any compatible topology. This was the original motivation 
for the development of a whole new field within general topology, 
convergence spaces.

Then in 1982, Alfred Fr\"olicher and Andreas Kriegl presented 
independently the solution to the question for the right differential 
calculus in infinite dimensions. They joined forces in the further 
development of the theory and the (up to now) final outcome is the 
book \cite{Fr\"olicher-Kriegl, 1988}.

In this section we will sketch the basic definitions and the most 
important results of the Fr\"olicher-Kriegl calculus.

\subheading{\nmb.{1.2}. The $c^\infty$-topology} Let $E$ be a 
locally convex vector space. A curve $c:\Bbb R\to E$ is called 
{\it smooth} or $C^\infty$ if all derivatives exist and are 
continuous - this is a concept without problems. Let 
$C^\infty(\Bbb R,E)$ be the space of smooth functions. It can be 
shown that $C^\infty(\Bbb R,E)$ does not depend on the locally convex 
topology of $E$, only on its associated bornology (system of bounded 
sets).

The final topologies with respect to the following sets of mappings 
into E coincide:
\roster
\item $C^\infty(\Bbb R,E)$.
\item Lipschitz curves (so that $\{\frac{c(t)-c(s)}{t-s}:t\neq s\}$ 
     is bounded in $E$). 
\item $\{E_B\to E: B\text{ bounded absolutely convex in }E\}$, where 
     $E_B$ is the linear span of $B$ equipped with the Minkowski 
     functional $p_B(x):= \inf\{\la>0:x\in\la B\}$.
\item Mackey-convergent sequences $x_n\to x$ (there exists a sequence 
     $0<\la_n\nearrow\infty$ with $\la_n(x_n-x)$ bounded).
\endroster
This topology is called the $c^\infty$-topology on $E$ and we write 
$c^\infty E$ for the resulting topological space. In general 
(on the space $\Cal D$ of test functions for example) it is finer 
than the given locally convex topology, it is not a vector space 
topology, since scalar multiplication is no longer jointly 
continuous. The finest among all locally convex topologies on $E$ 
which are coarser than $c^\infty E$ is the bornologification of the 
given locally convex topology. If $E$ is a Fr\'echet space, then 
$c^\infty E = E$. 

\subheading{\nmb.{1.3}. Convenient vector spaces} Let $E$ be a 
locally convex vector space. $E$ is said to be a {\it convenient 
vector space} if one of the following equivalent
conditions is satisfied (called $c^\infty$-completeness):
\roster
\item Any Mackey-Cauchy-sequence (so that $(x_n-x_m)$ is Mackey 
     convergent to 0) converges. 
\item If $B$ is bounded closed absolutely convex, then $E_B$ is a 
     Banach space.
\item Any Lipschitz curve in $E$ is locally Riemann integrable.
\item For any $c_1\in C^\infty(\Bbb R,E)$ there is 
     $c_2\in C^\infty(\Bbb R,E)$ with $c_1=c_2'$ (existence of 
     antiderivative).
\endroster

\proclaim{\nmb.{1.4}. Lemma} Let $E$ be a locally convex space.
Then the following properties are equivalent:
\roster
\item $E$ is $c^\infty$-complete.
\item If $f:\Bbb R\to E$ is scalarwise $\Lip^k$, then $f$ is 
     $\Lip^k$, for $k>1$.
\item If $f:\Bbb R\to E$ is scalarwise $C^\infty$ then $f$ is 
     differentiable at 0.
\item If $f:\Bbb R\to E$ is scalarwise $C^\infty$ then $f$ is 
     $C^\infty$.
\endroster
\endproclaim
Here a mapping $f:\Bbb R\to E$ is called $\Lip^k$ if all partial 
derivatives up to order $k$ exist and are Lipschitz, locally on 
$\Bbb R$. $f$ scalarwise $C^\infty$ means that $\la\o f$ is $C^\infty$  
for all continuous linear functionals on $E$.

This lemma says that on a convenient vector space one can recognize 
smooth curves by investigating compositions with continuous linear 
functionals.

\subheading{\nmb.{1.5}. Smooth mappings} Let $E$ and $F$ be locally 
convex vector spaces. A mapping $f:E\to F$ is called {\it smooth} or 
$C^\infty$, if $f\o c\in C^\infty(\Bbb R,F)$ for all 
$c\in C^\infty(\Bbb R,E)$; so 
$f_*: C^\infty(\Bbb R,E)\to C^\infty(\Bbb R,F)$ makes sense.
Let $C^\infty(E,F)$ denote the space of all smooth mapping from $E$ 
to $F$.

For $E$ and $F$ finite dimensional this gives the usual notion of 
smooth mappings: this has been first proved in \cite{Boman, 1967}.
Constant mappings are smooth. Multilinear mappings are smooth if and 
only if they are bounded. Therefore we denote by $L(E,F)$ the space 
of all bounded linear mappings from $E$ to $F$.

\subheading{\nmb.{1.6}. Structure on $C^\infty(E,F)$} We equip the 
space $C^\infty(\Bbb R,E)$ with the bornologification of the topology 
of uniform convergence on compact sets, in all derivatives 
separately. Then we equip the space $C^\infty(E,F)$ with the 
bornologification of the initial topology with respect to all 
mappings $c^*:C^\infty(E,F)\to C^\infty(\Bbb R,F)$, $c^*(f):=f\o c$, 
for all $c\in C^\infty(\Bbb R,E)$.

\proclaim{\nmb.{1.7}. Lemma } For locally convex spaces $E$ and $F$ 
we have:
\roster
\item If $F$ is convenient, then also $C^\infty(E,F)$ is convenient, 
     for any $E$. The space $L(E,F)$ is a closed linear subspace of 
     $C^\infty(E,F)$, so it also is convenient.
\item If $E$ is convenient, then a curve $c:\Bbb R\to L(E,F)$ is 
     smooth if and only if $t\mapsto c(t)(x)$ is a smooth curve in $F$ 
     for all $x\in E$.
\endroster
\endproclaim

\proclaim{\nmb.{1.8}. Theorem} The category of convenient vector 
spaces and smooth mappings is cartesian closed. So we have a natural 
bijection 
$$C^\infty(E\x F,G)\cong C^\infty(E,C^\infty(F,G)),$$
which is even a diffeomorphism.
\endproclaim

Of course this statement is also true for $c^\infty$-open subsets of 
convenient vector spaces. 

\proclaim{\nmb.{1.9}. Corollary } Let all spaces be convenient vector 
spaces. Then the following canonical mappings are smooth.
$$\align
&\operatorname{ev}: C^\infty(E,F)\x E\to F,\quad 
     \operatorname{ev}(f,x) = f(x)\\
&\operatorname{ins}: E\to C^\infty(F,E\x F),\quad
     \operatorname{ins}(x)(y) = (x,y)\\
&(\quad)^\wedge :C^\infty(E,C^\infty(F,G))\to C^\infty(E\x F,G)\\
&(\quad)\spcheck :C^\infty(E\x F,G)\to C^\infty(E,C^\infty(F,G))\\
&\operatorname{comp}:C^\infty(F,G)\x C^\infty(E,F)\to C^\infty(E,G)\\
&C^\infty(\quad,\quad):C^\infty(F,F')\x C^\infty(E',E)\to 
     C^\infty(C^\infty(E,F),C^\infty(E',F'))\\
&\qquad (f,g)\mapsto(h\mapsto f\o h\o g)\\
&\prod:\prod C^\infty(E_i,F_i)\to C^\infty(\prod E_i,\prod F_i)
\endalign$$
\endproclaim

\proclaim{\nmb.{1.10}. Theorem} Let $E$ and $F$ be convenient vector 
spaces. Then the differential operator 
$$\gather d: C^\infty(E,F)\to C^\infty(E,L(E,F)), \\
df(x)v:=\lim_{t\to0}\frac{f(x+tv)-f(x)}t,
\endgather$$
exists and is linear and bounded (smooth). Also the chain rule holds: 
$$d(f\o g)(x)v = df(g(x))dg(x)v.$$
\endproclaim

\subheading{\nmb.{1.11}. Remarks } Note that the conclusion of 
theorem \nmb!{1.8} is the starting point of the classical calculus of 
variations, where a smooth curve in a space of functions was assumed 
to be just a smooth function in one variable more.

If one wants theorem \nmb!{1.8} to be true and assumes some other obvious 
properties, then the calculus of smooth functions is already uniquely 
determined.

There are, however, smooth mappings which are not continuous. This is 
unavoidable and not so horrible as it might appear at first sight. 
For example the evaluation $E\x E'\to\Bbb R$ is jointly continuous if 
and only if $E$ is normable, but it is always smooth. Clearly smooth 
mappings are continuous for the $c^\infty$-topology.

For Fr\'echet spaces smoothness in the sense described here coincides 
with the notion $C^\infty_c$ of \cite{Keller, 1974}. This is the 
differential calculus used by \cite{Michor, 1980}, \cite{Milnor, 
1984}, and \cite{Pressley-Segal, 1986}.

A prevalent opinion in contemporary mathematics is,
that for infinite
dimensional calculus each serious application needs its own 
foundation. By a serious application one obviously means some 
application of
a hard inverse function theorem. These theorems can be proved,
if by assuming enough a priori estimates one creates enough
Banach space situation for some modified iteration procedure to
converge. Many authors try to build their platonic idea of an a
priori estimate into their differential calculus. We think that this 
makes the calculus
inapplicable and hides the origin of the a priori estimates. We
believe, that the calculus itself should be as easy to use as
possible, and that all further assumptions (which most often
come from ellipticity of some nonlinear partial differential
equation of geometric origin) should be treated separately, in a
setting depending on the specific problem. We are sure that in
this sense the Fr\"olicher-Kriegl calculus as presented here and its 
holomorphic and real analytic offsprings in sections \nmb!{2} and 
\nmb!{3} below are universally usable for most applications.

We believe that the recent development of the theory of locally 
convex spaces missed its original aim: the development of calculus.
It laid too much emphasis on (locally convex) topologies and ignored 
and denigraded
the work of Hogbe-Nlend, Colombeau and collaborators on the original idea of 
Sebasti\~ao e Silva, that bornologies are the right concept for this 
kind of functional analysis.

\heading\nmb0{2}. Calculus of holomorphic mappings \endheading

\subheading{\nmb.{2.1}} Along the lines of thought of the 
Fr\"olicher-Kriegl calculus of smooth mappings, in \cite{Kriegl-Nel, 
1985} the cartesian closed setting for holomorphic mappings was 
developed. The right definition of this calculus was already given 
by \cite{Fantappi\'e, 1930 and 1933}. We will now sketch the basics 
and the main results. It can be shown that again convenient vector 
spaces are the right ones to consider. Here we will start with them 
for the sake of shortness.

\subheading{\nmb.{2.2}} Let $E$ be a complex locally convex vector 
space whose underlying real space is convenient -- this will be 
called convenient in the sequel. Let $\Bbb D\subset \Bbb C$ be the 
open unit disk and let us denote by $\Cal H(\Bbb D,E)$ the space of 
all mappings $c:\Bbb D\to E$ such that $\la\o c:\Bbb D\to \Bbb C$ is 
holomorphic for each continuous complex-linear functional $\la$ on 
$E$. Its elements will be called the holomorphic curves.

If $E$ and $F$ are convenient complex vector spaces (or 
$c^\infty$-open sets therein), a mapping 
$f:E\to F$ is called {\it holomorphic} if $f\o c$ is a holomorphic 
curve in $F$ for each holomorphic curve $c$ in $E$. Obviously $f$ is 
holomorphic if and only if $\la\o f:E\to \Bbb C$ is holomorphic for 
each complex linear continuous functional $\la$ on $F$. Let 
$\Cal H(E,F)$ denote the space of all holomorphic mappings from $E$ to 
$F$. 

\proclaim{\nmb.{2.3}. Theorem (Hartog's theorem)} Let $E_k$ for 
$k=1,2$ and $F$ be complex convenient vector spaces and let 
$U_k\subset E_k$ be $c^\infty$-open. A mapping $f:U_1\x U_2\to F$ is 
holomorphic if and only if it is separately holomorphic (i\. e\. 
$f(\quad,y)$ and $f(x,\quad)$ are holomorphic for all $x\in U_1$ and 
$y\in U_2$).
\endproclaim

This implies also that in finite dimensions we have recovered the 
usual definition.

\proclaim{\nmb.{2.4} Lemma} If $f:E\supset U\to F$ is holomorphic 
then $df:U\x E\to F$ exists, is holomorphic and $\Bbb C$-linear in 
the second variable. 

A multilinear mapping is holomorphic if and only if it is bounded.
\endproclaim

\proclaim{\nmb.{2.5} Lemma} If $E$ and $F$ are Banach spaces and $U$ 
is open in $E$, then for a mapping $f:U\to F$ the following 
conditions are equivalent:
\roster
\item $f$ is holomorphic.
\item $f$ is locally a convergent series of homogeneous continuous 
     polynomials.
\item $f$ is $\Bbb C$-differentiable in the sense of Fr\'echet.
\endroster 
\endproclaim

\proclaim{\nmb.{2.6} Lemma} Let $E$ and $F$ be convenient vector 
spaces. A mapping $f:E\to F$ is holomorphic if and only if it is 
smooth and its derivative is everywhere $\Bbb C$-linear.
\endproclaim

An immediate consequence of this result is that $\Cal H(E,F)$ is a 
closed linear subspace of $C^\infty(E_{\Bbb R},F_{\Bbb R})$ and so it 
is a convenient vector space if $F$ is one, by \nmb!{1.7}. The chain 
rule follows from \nmb!{1.10}. The following theorem is an easy 
consequence of \nmb!{1.8}.

\proclaim{\nmb.{2.7} Theorem} The category of convenient complex 
vector spaces and holomorphic mappings between them is cartesian 
closed, i\. e\.
$$\Cal H(E\x F,G) \cong \Cal H(E,\Cal H(F,G)).$$
\endproclaim

An immediate consequence of this is again that all canonical 
structural mappings as in \nmb!{1.9} are holomorphic.

\heading\nmb0{3}. Calculus of real analytic mappings \endheading 

\subheading{\nmb.{3.1}} In this section we sketch the cartesian closed 
setting to real analytic mappings in infinite dimension following the 
lines of the Fr\"olicher-Kriegl calculus, as it is presented in 
\cite{Kriegl-Michor, 1990}. Surprisingly enough one has to deviate 
from the most obvious notion of real analytic curves in order to get 
a meaningful theory, but again convenient vector spaces turn out to 
be the right kind of spaces.

\subheading{\nmb.{3.2}. Real analytic curves} Let $E$ be a real 
convenient vector space with dual $E'$. A curve $c:\Bbb R\to E$ is 
called {\it real analytic} if $\la\o c:\Bbb R\to \Bbb R$ is real 
analytic for each $\la\in E'$.
It turns out that the set of these curves depends only on the 
bornology of $E$.

In contrast a curve is called {\it topologically real analytic} if it 
is locally given by power series which converge in the topology of 
$E$. They can be extended to germs of holomorphic curves along $\Bbb R$ 
in the complexification $E_{\Bbb C}$ of $E$. If the dual $E'$ of $E$ 
admits a Baire topology which is compatible with the duality, then 
each real analytic curve in $E$ is in fact topologically real analytic 
for the bornological topology on $E$.

\subheading{\nmb.{3.3}. Real analytic mappings} Let $E$ and $F$ be 
convenient vector spaces. Let $U$ be a $c^\infty$-open set in $E$. A 
mapping $f:U\to F$ is called {\it real analytic} if and only if it is 
smooth (maps smooth curves to smooth curves) and maps real analytic 
curves to real analytic curves. 

Let $C^\omega(U,F)$ denote the space of all real analytic mappings. 
We equip the space
$C^\om(U,\Bbb R)$ of all real analytic functions 
with the initial topology
with respect to the families of mappings
$$\gather C^\om(U,\Bbb R) @>{c^*}>> C^\om(\Bbb R,\Bbb R),\text{
for all }c\in C^\om(\Bbb R,U)\\ 
C^\om(U,\Bbb R) @>{c^*}>> C^\infty(\Bbb R,\Bbb R),\text{
for all }c\in C^\infty(\Bbb R,U),
\endgather$$
where $C^\infty(\Bbb R,\Bbb R)$ carries the topology of compact 
convergence in each derivative separately as in section \nmb!{1}, and 
where $C^\omega(\Bbb R,\Bbb R)$ is equipped with the final locally 
convex topology 
with respect to the embeddings (restriction mappings) of all spaces 
of holomorphic mappings from a neighborhood $V$ of $\Bbb R$ in 
$\Bbb C$ mapping $\Bbb R$ to $\Bbb R$, and each of these spaces 
carries the topology of compact convergence.

Furthermore we equip the space
$C^\om(U,F)$ with the initial topology with respect to
the family of mappings 
$$ C^\om(U,F) @>{\la_*}>> C^\om(U,\Bbb R),\text{ for all }\la\in
F'.$$
It turns out that this is again a convenient space.

\proclaim{\nmb.{3.4}. Theorem} In the setting of \nmb!{3.3} a mapping 
$f:U\to F$ is real analytic if and only if it is smooth and is real 
analytic along each affine line in $E$.
\endproclaim

\proclaim{\nmb.{3.5}. Theorem} The category of convenient spaces and 
real analytic mappings is cartesian closed. So the equation
$$C^\omega(U,C^\omega(V,F))\cong C^\omega(U\x V,F)$$
is valid for all $c^\infty$-open sets $U$ in $E$ and $V$ in $F$, 
where $E$, $F$, and $G$ are convenient vector spaces.
\endproclaim

This implies again that all structure mappings as in \nmb!{1.9} are 
real analytic. Furthermore the differential operator 
$$d:C^\omega(U,F)\to C^\omega(U,L(E,F))$$ 
exists, is unique and real 
analytic. Multilinear mappings are real analytic if and only if they 
are bounded. Powerful real analytic uniform boundedness principles 
are available.

\heading\nmb0{4}. Infinite dimensional manifolds \endheading

\subheading{\nmb.{4.1}} In this section we will concentrate on two 
topics: Smooth partitions of unity, and several kinds of tangent 
vectors.

\subheading{\nmb.{4.2}} In the usual way we define manifolds by 
gluing $c^\infty$-open sets in convenient vector spaces via smooth 
(holomorphic, real analytic) diffeomorphisms. Then we equip them with 
the identification topology with respect to the $c^\infty$-topologies 
on all modeling spaces. We require some properties from this 
topology, like Hausdorff and regular (which here is not a consequence 
of Hausdorff).

Mappings between manifolds are smooth (holomorphic, real analytic), 
if they have this property when composed which any chart mappings.

\proclaim{\nmb.{4.3}. Lemma} A manifold $M$ is metrizable if
and only if it is paracompact and modeled on Fr\'echet spaces.
\endproclaim

\proclaim{\nmb.{4.4}. Lemma}
For a convenient vector space $E$ the
set $C^\infty(M,E)$ of smooth $E$-valued functions on a manifold
$M$ is again a convenient vector space. Likewise for the real
analytic and holomorphic case.
\endproclaim

\proclaim{\nmb.{4.5}. Theorem}
If $M$ is a smooth manifold modeled on convenient vector spaces
admitting smooth bump functions and
$\Cal U$ is a locally finite open cover of $M$, then there
exists a smooth partition of unity
$\{\ph_U:U\in\Cal U\}$ with $carr/supp(\ph_U)\subset U$ for all
$U\in\Cal U$.

If $M$ is in addition paracompact, then this is true for every
open cover $\Cal U$ of $M$.

Convenient vector spaces which are nuclear admit smooth bump 
functions.
\endproclaim

\subheading{\nmb.{4.6}. The tangent spaces of a convenient vector space $E$} 
Let $a \in E$. A {\it kinematic tangent vector} with foot
point $a$ is simply 
a pair $(a,X)$ with $X\in E$. Let $T_aE=E$ be the space of all
kinematic tangent vectors with foot point $a$. It consists of
all derivatives $c'(0)$ at 0 of smooth curves $c:\Bbb R\to E$
with $c(0)=a$, which explains the choice of the name kinematic.

For each open neighborhood $U$ of $a$ in $E$
$(a,X)$ induces a linear mapping $X_a:C^\infty(U,\Bbb R) \to  \Bbb R$
by $X_a(f) := df(a)(X_a)$, which is continuous for the convenient
vector space topology on $C^\infty(U,\Bbb R)$, and satisfies 
$X_a(f\cdot g) = X_a(f)\cdot g(a) + f(a)\cdot X_a(g)$,
so it is a {\it continuous derivation over $ev_a$}. 
The value $X_a(f)$ depends only on the germ of $f$ at $a$.

An {\it operational tangent vector} of $E$ with foot
point $a$ is a bounded derivation 
$\partial :C^\infty_a(E,\Bbb R)\to \Bbb R$ over $ev_a$. 
Let $D_aE$ be the vector space of all these derivations.
Any $\partial\in D_aE$ induces a bounded derivation
$C^\infty(U,\Bbb R)\to \Bbb R$ over $ev_a$ for each open
neighborhood $U$ of $a$ in $E$. So the vector space
$D_aE$ is a closed linear subspace of the convenient vector
space $\prod_U L(C^\infty(U,\Bbb R),\Bbb R)$. We equip $D_aE$
with the induced convenient vector space structure. Note that
the spaces $D_aE$ are isomorphic for all $a\in E$.

\subheading{Example} Let $Y\in E''$ be an element in the bidual
of $E$. Then for each $a\in E$ we have an operational tangent vector
$Y_a\in D_aE$, given by $Y_a(f):= Y(df(a))$. So we have a
canonical injection $E''\to D_aE$.
 
\subheading{Example} Let $\ell: L^2(E;\Bbb R)\to \Bbb R$ be a bounded
linear functional which vanishes on the subset $E'\otimes E'$.
Then for each $a\in E$ we have an operational tangent vector 
$\partial^2_{\ell}|_a\in D_aE$ given by
$\partial^2_\ell|_a(f) := \ell(d^2f(a))$, since
$$\align &\ell(d^2(fg)(a))=\\
&\qquad = \ell(d^2f(a)g(a)+df(a)\otimes
    dg(a)+dg(a)\otimes df(a)+ f(a)d^2g(a))\\
&\qquad = \ell(d^2f(a))g(a) + 0 + f(a)\ell(d^2g(a)).
\endalign$$

\proclaim{\nmb.{4.7}. Lemma} Let $\ell\in L^k_{sym}(E;\Bbb
R)'$ be a bounded linear functional which vanishes on the subspace
$$\sum_{i=1}^{k-1}L^i_{sym}(E;\Bbb R)\vee L^{k-i}_{sym}(E;\Bbb R)$$
of decomposable elements of $L^k_{sym}(E;\Bbb R)$. Then $\ell$
defines an operational tangent vector $\partial^k_\ell|_a\in
D_aE$ for each $a\in E$ by
$$\partial^k_\ell|_a(f) := \ell(d^kf(a)).$$
The linear mapping $\ell\mapsto \partial^k_\ell|_a$  is an
embedding onto a topological direct summand $D^{(k)}_aE$ of $D_aE$.
Its left inverse is given by 
$\partial\mapsto (\Ph\mapsto \partial((\Ph\o diag)(a+\quad)))$.
The  sum $\sum_{k>0}D^{(k)}_aE$ in $D_aE$ is a direct one.  
\endproclaim

\proclaim{\nmb.{4.8}. Lemma} If $E$ is an infinite
dimensional Hilbert space, all operational tangent space summands
$D^{(k)}_0E$ are not zero.
\endproclaim

\subheading{\nmb.{4.9}. Definition} A convenient vector space
is said to have the {\it (borno\-logical) approximation
property} if $E'\otimes E$ is dense in $L(E,E)$ in the
bornological locally convex topology.

The following spaces have the bornological approximation property:
$\Bbb R^{(\Bbb N)}$, nuclear Fr\'echet spaces, nuclear (LF) spaces.

\proclaim{\nmb.{4.10} Theorem} Let $E$ be a convenient vector
space which has the approximation property. 
Then we have $D_aE=D^{(1)}_aE\cong E''$. So if $E$ is in
addition reflexive, each operational tangent vector is a
kinematic one.
\endproclaim 

\subheading{\nmb.{4.11}} The kinematic tangent bundle $TM$ of a 
manifold $M$ is constructed by gluing all the kinematic tangent 
bundles of charts with the help of the kinematic tangent mappings 
(derivatives) of 
the chart changes. $TM\to M$ is a vector bundle and 
$T:C^\infty(M,N)\to C^\infty(TM,TN)$ is well defined and has the 
usual properties.

\subheading{\nmb.{4.12}} The operational tangent bundle $DM$ of a manifold
$M$ is constructed by gluing all operational tangent spaces of 
charts.
Then $\pi_M:DM \to M$ is again a vector bundle which contains 
the kinematic tangent bundle $TM$ as a splitting
subbundle. Also for each $k\in\Bbb N$ the same gluing
construction as above gives us tangent bundles $D^{(k)}M$ which
are splitting sub bundles of $DM$. The mappings
$D^{(k)}:C^\infty(M,N)\to C^\infty(D^{(k)}M,D^{(k)}N)$ are well 
defined for all $k$ (including no $k$) and have the usual properties.

Note that for manifolds modeled on reflexive spaces having the 
borno\-logical approximation property the operational and the kinematic 
tangent bundles coincide.

\heading\nmb0{5}. Manifolds of mappings \endheading

\proclaim{\nmb.{5.1}. Theorem (Manifold structure of
$C^\infty(M,N)$)} Let $M$ and $N$ be smooth finite dimensional
manifolds, let $M$ be compact. Then the space $C^\infty(M,N)$ of
all smooth mappings from $M$ to $N$ is a smooth
manifold, modeled on spaces $C^\infty(f^*TN)$ of smooth
sections of pullback bundles along $f:M\to N$ over $M$.
\endproclaim
A careful description of this theorem (but without the help of the 
Fr\"olicher-Kriegl calculus) can be found in \cite{Michor, 1980}.
We include a proof of this result here because the result is important 
and the proof is much simpler now.
\demo{Proof} Choose a smooth Riemannian metric on $N$.
Let $\exp: TN\supseteq U\to N$ 
be the smooth exponential mapping of this Riemannian metric,
defined on a suitable open neighborhood of the zero section. We
may assume that $U$ is chosen in such a way that
$(\pi_N,\exp):U\to N\x N$ is a smooth diffeomorphism onto
an open neighborhood $V$ of the diagonal.

For $f\in C^\infty(M,N)$ we consider the pullback vector bundle
$$\CD M\x_NTN @= f^*TN     @>{\pi_N^*f}>>           TN \\
@.               @V{f^*\pi_N}VV                @VV{\pi_N}V\\
              @. M         @>f>>              N.  \endCD$$
Then $C^\infty(f^*TN)$ is canonically isomorphic to 
$C^\infty_f(M,TN):= \{h\in C^\infty(M,TN):\pi_N\o h=f\}$ via 
$s\mapsto (\pi_N^*f)\o s$ and $(Id_M,h)\gets h$. We consider the 
space $C^\infty_c(f^*TN)$ of all smooth sections with compact support 
and equip it with the inductive limit topology
$$C^\infty_c(f^*TN)=\injlim_K C^\infty_K(f^*TN),$$
where $K$ runs through all compact sets in $M$ and each of the spaces
$C^\infty_K(f^*TN)$ is equipped with the topology of uniform 
convergence (on $K$) in all derivatives separately.
Now let 
$$\gather 
U_f :=\{g\in C^\infty(M,N):(f(x),g(x))\in V
     \text{ for all }x\in M,g\sim f\},\\
u_f:U_f\to C^\infty_c(f^*TN),\\
u_f(g)(x) = (x,\exp_{f(x)}\i(g(x))) = (x,((\pi_N,\exp)\i\o(f,g))(x)).
\endgather$$
Here $g\sim f$ means that $g$ equals $f$ off some compact set.
Then $u_f$ is a bijective mapping from $U_f$ onto the set
$\{s\in C^\infty_c(f^*TN): s(M)\subseteq f^*U\}$, whose inverse is
given by $u_f\i(s) = \exp\o(\pi_N^*f)\o s$, where we view 
$U \to N$ as fiber bundle. The set $u_f(U_f)$ is open in
$C^\infty_c(f^*TN)$ for the topology described above.

Now we consider the atlas $(U_f,u_f)_{f\in C^\infty(M,N)}$ for
$C^\infty(M,N)$. Its chart change mappings are given 
for $s\in u_g(U_f\cap U_g)\subseteq C^\infty_c(g^*TN)$ by 
$$\align 
(u_f\o u_g\i)(s) &= (Id_M,(\pi_N,\exp)\i\o(f,\exp\o(\pi_N^*g)\o s)) \\
&= (\tau_f\i\o\tau_g)_*(s),
\endalign$$ 
where $\tau_g(x,Y_{g(x)}) := (x,\exp_{g(x)}(Y_{g(x)})))$
is a smooth diffeomorphism 
$\tau_g:g^*TN \supseteq g^*U \to (g\x Id_N)\i(V)\subseteq M\x N$
which is fiber respecting over $M$. 

Smooth curves in $C^\infty_c(f^*TN)$ are just smooth sections of the bundle
$pr_2^*f^*TN\to \Bbb R\x M$, which have compact support in $M$ 
locally in$\Bbb R$. 
The chart change
$u_f\o u_g\i = (\tau_f\i\o \tau_g)_*$ is defined on an open
subset and obviously maps smooth curves to smooth curves, therefore 
it is also smooth.

Finally we put the identification topology from this atlas onto
the space $C^\infty(M,N)$, which is obviously finer than the
compact open topology and thus Hausdorff.

The equation $u_f\o u_g\i = (\tau_f\i\o \tau_g)_*$ shows that
the smooth structure does not depend on the choice of the
smooth Riemannian metric on $N$.
\qed\enddemo

\proclaim{\nmb.{5.2}. Theorem ($C^\om$-manifold structure of
$C^\om(M,N)$)} Let $M$ and $N$ be real analytic 
manifolds, let $M$ be compact. Then the space $C^\om(M,N)$ of
all real analytic mappings from $M$ to $N$ is a real analytic
manifold, modeled on spaces $C^\om(f^*TN)$ of real analytic
sections of pullback bundles along $f:M\to N$ over $M$.
\endproclaim
The proof can be found in \cite{Kriegl-Michor, 1990}. It is a variant 
of the above proof, using a real analytic Riemannian metric.

\proclaim{\nmb.{5.3}. Theorem ($C^\om$-manifold structure on
$C^\infty(M,N)$)} Let $M$ and $N$ be real analytic
manifolds, with $M$ compact. Then the smooth manifold
$C^\infty(M,N)$ is even a real analytic manifold. 
\endproclaim
\demo{Proof} For a fixed real analytic exponential mapping on $N$
the charts $(U_f,u_f)$ from \nmb!{5.1} 
for $f\in C^\om(M,N)$ form a smooth atlas
for $C^\infty(M,N)$, since $C^\om(M,N)$ is dense in
$C^\infty(M,N)$ by \cite {Grauert, 1958, Proposition 8}. The chart
changings $u_f\o u_g\i = (\tau_f\i\o \tau_g)_*$ are real
analytic: this follows from a careful description of the set of real 
analytic curves into $C^\infty(f^*TN)$. See again 
\cite{Kriegl-Michor, 1990, 7.7} for more details.
\qed\enddemo

\subheading{\nmb.{5.4} Remark} If $M$ is not compact,
$C^\om(M,N)$ is dense in $C^\infty(M,N)$ for the
Whitney-$C^\infty$-topology by \cite{Grauert, 1958, Prop. 8}.
This is not the case for the topology used in \nmb!{5.1} in which
$C^\infty(M,N)$ is a smooth manifold. The charts $U_f$ for $f\in
C^\om(M,N)$ do not cover $C^\infty(M,N)$.

\proclaim{\nmb.{5.5}. Theorem} Let $M$
and $N$ be smooth manifolds. Then the
two infinite dimensional smooth 
vector bundles $TC^\infty(M,N)$ and $C^\infty(M,TN)$ over $C^\infty(M,N)$
are canonically isomorphic. The same assertion is true for
$C^\om(M,N)$, if $M$ is compact. 
\endproclaim

\proclaim{\nmb.{5.6}. Theorem (Exponential law)} Let $\Cal M$ be
a (possibly infinite 
dimensional) smooth manifold, and let $M$ and $N$ be
finite dimensional smooth manifolds.

Then we  have a canonical embedding 
$$C^\infty(\Cal M, C^\infty(M,N))\subseteq C^\infty(\Cal M\x M,N),$$ 
where we have equality if and only if $M$ is compact.

If $M$ and $N$ are real analytic manifolds with $M$ compact we have
$$C^\om(\Cal M,C^\omega(M,N))= C^\omega(\Cal M\x M,N)$$
for each real analytic (possibly infinite dimensional) manifold.
\endproclaim

\proclaim{\nmb.{5.7}. Corollary} If $M$ is compact and $M$, $N$ are
finite dimensional smooth manifolds, then the evaluation
mapping $\ev: C^\infty(M,N)\x M\to N$ is smooth.

If $P$ is another compact smooth manifold, then the
composition mapping
$\operatorname{comp}:C^\infty(M,N)\x C^\infty(P,M)\to C^\infty(P,N)$ 
is smooth.

In particular $f_*:C^\infty(M,N)\to C^\infty(M,N')$ and 
$g^*:C^\infty(M,N)\to C^\infty(P,N)$ are smooth for 
$f\in C^\infty(N,N')$ and $g\in C^\infty(P,M)$.

The corresponding statement for real analytic mappings is also true.
\endproclaim

\proclaim{\nmb.{5.8}. Theorem (Diffeomorphism groups)} 
For a smooth manifold 
$M$ the group $\operatorname{Diff}(M)$ of all smooth
diffeomorphisms of $M$ is an open submani\-fold of $C^\infty(M,M)$,
composition and inversion are smooth.

The Lie algebra of the smooth infinite dimensional
Lie group $\operatorname{Diff}(M)$ is the convenient vector
space $C^\infty_c(TM)$ of all smooth vector fields on $M$ with compact 
support,
equipped with the negative of the usual Lie bracket. The
exponential mapping $\operatorname{Exp}:C^\infty_c(TM)\to
\operatorname{Diff}^\infty(M)$ is the flow mapping to time 1, and
it is smooth.

For a compact real analytic manifold $M$ the group 
$\operatorname{Diff}^\om(M)$ of all real analytic diffeomorphisms is 
a real analytic Lie group with Lie algebra $C^\om(TM)$ and with real 
analytic exponential mapping.
\endproclaim

\subheading{\nmb.{5.9}. Remarks}
The group $\operatorname{Diff}(M)$ of smooth diffeomorphisms does not
carry any real analytic Lie group structure by
\cite{Milnor, 1984, 9.2}, and it has no complexification in
general, see \cite{Pressley-Segal, 1986, 3.3}.  The mapping
$$Ad\o\operatorname{Exp}:C^\infty_c(TM)\to \operatorname{Diff}(M)\to
L(C^\infty(TM),C^\infty(TM))$$ 
is not real analytic, see \cite{Michor, 1983, 4.11}.

For $x\in M$ the mapping 
$ev_x\o \operatorname{Exp}:C^\infty_c(TM)\to
\operatorname{Diff}(M)\to M$ is not real analytic since 
$(ev_x\o\operatorname{Exp})(tX) = Fl^X_t(x)$, which is not real
analytic in $t$ for general smooth $X$.

The exponential mapping $\operatorname{Exp}:C^\infty_c(TM)\to
\operatorname{Diff}(M)$ is in a very strong
sense not surjective: In \cite{Grabowski, 1988} it is shown, that
$\operatorname{Diff}(M)$ contains an arcwise connected free
subgroup on $2^{\aleph_0}$ generators which meets the image of
$\operatorname{Exp}$ only at the identity. 

The real analytic Lie group $\operatorname{Diff}^\om(M)$ is {\it
regular} in the sense of \cite{Milnor, 1984. 7.6}, who weakened
the original concept of \cite{Omori, 1982}. This condition means
that the mapping associating the evolution operator to each time
dependent vector field on $M$ is smooth. It is even real
analytic, compare the proof of theorem \nmb!{5.9}.

\proclaim{\nmb.{5.10}. Theorem}
Let $M$ and $N$ be 
smooth manifolds. Then the diffeomorphism group 
$\operatorname{Diff}(M)$ acts smoothly from the right on the smooth
manifold $\operatorname{Imm}(M,N)$
of all smooth immersions $M\to N$, which is an open subset
of $C^\infty(M,N)$. 

Then the space of orbits 
$\operatorname{Imm}(M,N)/\operatorname{Diff}(M)$
is Hausdorff in the quotient topology.

Let $\operatorname{Imm}_{\text{free}}(M,N)$ be set of all immersions, 
on which $\operatorname{Diff}(M)$ acts freely. Then this is open in 
$C^\infty(M,N)$ and is the total space of a smooth principal fiber 
bundle 
$$\operatorname{Imm}_{\text{free}}(M,N)\to 
\operatorname{Imm}_{\text{free}}(M,N)/\operatorname{Diff}(M).$$ 

In particular the space $\operatorname{Emb}(M,N)$ of all smooth 
embeddings is the total space of smooth principal fiber bundle.
\endproclaim

This is proved in \cite{Cervera-Mascaro-Michor, 1989}, where also the 
existence of smooth transversals to each orbit is shown and the 
stratification of the orbit space into smooth manifolds is given.

\proclaim{\nmb.{5.11}. Theorem (Principal bundle of embeddings)}
Let $M$ and $N$ be real 
analytic manifolds with $M$ compact. Then the set $Emb^\om(M,N)$
of all real analytic embeddings $M\to N$ is an open submanifold
of $C^\om(M,N)$. It is the total space of a real analytic
principal fiber bundle with structure group
$\operatorname{Diff}^\om(M)$, whose real analytic base manifold
is the space of all submanifolds of $N$ of type $M$.
\endproclaim

See \cite{Kriegl-Michor, 1990}, section 6.

\proclaim{\nmb.{5.12}. Theorem (Classifying space for
$\operatorname{Diff}(M)$)}  Let $M$ be a compact smooth
manifold. Then the space $Emb(M,\ell^2)$ of
smooth embeddings of $M$ into the
Hilbert space $\ell^2$ is the total space of a smooth principal fiber
bundle with structure group $\operatorname{Diff}(M)$ and
smooth base manifold $B(M,\ell^2)$, which is a classifying
space for the Lie group $\operatorname{Diff}(M)$.
It carries a universal $\operatorname{Diff}(M)$-connection.

In other words: 
$$(\operatorname{Emb}(M,\ell^2)\x_{\operatorname{Diff}(M)}M\to B(M,\ell^2)$$ 
classifies fiber bundles with typical fiber $M$ and carries a 
universal (generalized) connection.
\endproclaim

See \cite{Michor, 1988, section 6}.

\heading\nmb0{6}. Manifolds for algebraic topology \endheading

\subheading{\nmb.{6.1} Convention} In this section the space
$\Bbb R^{(\Bbb N)}$ of all finite sequences 
with the direct sum topology will be denoted by $\Bbb R^\infty$
following the common usage in algebraic topology. It is a
convenient vector space. 

We consider on it the weak inner product
$\langle x,y\rangle := \sum x_iy_i$, which is bilinear and
bounded, therefore smooth. It is called weak, since it is non
degenerate in the following sense: the
associated linear mapping 
$\Bbb R^{\infty}\to (\Bbb R^{\infty})'= \Bbb R^{\Bbb N}$ 
is injective, but far from being surjective. 
We will also use the weak Euclidean distance
$|x|:=\sqrt{\langle x,x\rangle}$, whose square is a smooth function. 

\subheading{\nmb.{6.2}. Example: The sphere $S^\infty$} 

The {\it sphere $S^\infty$} is the
set $\{x\in R^{\infty}: \langle x,x\rangle=1\}$. This is the
usual infinite dimensional sphere used in algebraic topology,
the topological inductive limit of $S^n\subset S^{n+1}\subset\ldots$.

We show that $S^\infty$ is a smooth manifold by
describing an
explicit smooth atlas for $S^{\infty}$, the {\it stereographic atlas}.
Choose $a\in S^{\infty}$ ("south pole"). Let 
$$\alignat3
&U_+:= S^\infty\setminus\{a\},&\qquad &u_+:U_+\to\{a\}^\bot, &\qquad
&u_+(x)=\tfrac{x-\langle x,a\rangle a}{1-\langle x,a\rangle},\\
&U_-:= S^\infty\setminus\{-a\},&\qquad &u_-:U_-\to\{a\}^\bot, &\qquad
&u_-(x)=\tfrac{x-\langle x,a\rangle a}{1+\langle x,a\rangle}.
\endalignat$$
 From an obvious drawing in the 2-plane through 0, $x$, and $a$
it is easily seen that $u_+$ is the usual stereographic
projection. We also get 
$$u_+\i(y) = \tfrac{|y|^2-1}{|y|^2+1}a+\tfrac 2{|y|^2+1}y
\qquad\text{for }y\in \{a\}^\bot\setminus\{0\}$$
and $(u_-\o u_+\i)(y) = \frac y{|y|^2}$. The latter equation can
directly be seen from the drawing using "Strahlensatz".

The two stereographic charts above can be extended to charts on
open sets in $\Bbb R^{\infty}$ in such a way that $S^\infty$
becomes a splitting submanifold of $\Bbb R^{\infty}$:
$$\gather 
\tilde u_+:\Bbb R^{\infty}\setminus [0,+\infty)a \to a^\bot + (-1,+\infty)a\\
\tilde u_+(z) := u_+(\frac z{|z|}) + (|z|-1)a.
\endgather$$

Since the model space $\Bbb R^{\infty}$ of $S^\infty$ has the
bornological approximation property by \nmb!{4.9}, and is
reflexive, by \nmb!{4.10} the operational tangent bundle of
$S^\infty$ equals the kinematic one: $DS^\infty = TS^\infty$.

We claim that $TS^\infty$ is diffeomorphic to 
$\{(x,v)\in S^\infty\x\Bbb R^{\infty}: \langle x,v\rangle=0\}$.

The $X_x\in T_xS^\infty$ are exactly of the form $c'(0)$ for a smooth
curve $c:\Bbb R\to S^\infty$ with $c(0)= x$ by \nmb!{4.11}. Then 
$0=\tfrac d{dt}|0\langle c(t),c(t)\rangle = 2\langle x,X_x\rangle$.
For $v\in x^\bot$ we use $c(t)= \cos(|v|t)x + \sin(|v|t)\tfrac v{|v|}$.

The construction of $S^\infty$ works for any positive definite bounded
bilinear form on any convenient vector space.

\subheading{\nmb.{6.3}. Example. The Grassmannians
and the Stiefel manifolds}\newline
The {\it Grassmann manifold
$G(k,\infty;\Bbb R)$} is the set of all
k-dimen\-sio\-nal linear subspaces of the space of all finite
sequences $\Bbb R^{\infty}$. 

The {\it Stiefel manifold
$O(k,\infty;\Bbb R)$ of orthonormal $k$-frames} is the set 
of all linear isometries $\Bbb R^k\to \Bbb R^{\infty}$, where
the latter space is again equipped with the standard weak inner
product described at the beginning of \nmb!{6.2}. 

The {\it Stiefel manifold $GL(k,\infty;\Bbb R)$ of all
$k$-frames} is the set of all
injective linear mappings $\Bbb R^k\to \Bbb R^{\infty}$.

There is a canonical {\it transposition mapping
$(\quad)^t:L(\Bbb R^k,\Bbb R^{\infty})\to L(\Bbb R^{\infty},\Bbb R^k)$} 
which is given by
$$A^t: \Bbb R^{\infty} @>incl>> \Bbb R^{\Bbb N}=
\left(\Bbb R^{\infty}\right)' @>A'>> (\Bbb R^k)'= \Bbb R^k $$
and satisfies 
$\langle A^t(x),y\rangle = \langle x,A(y)\rangle$.
The transposition mapping is boun\-ded and linear, so it is real analytic.

Then we have
$$GL(k,\infty)=\{A\in L(\Bbb R^k,\Bbb R^{\infty}): A^t\o A\in GL(k)\},$$
since $A^t\o A\in GL(k)$ if and only if 
$\langle Ax,Ay\rangle = \langle A^tAx,y\rangle=0$ for all $y$
implies $x=0$, which is equivalent to $A$ injective. So in
particular $GL(k,\infty)$ is open in $L(\Bbb R^k,\Bbb R^\infty)$.
The Lie group $GL(k)$ acts
freely from the right on the space $GL(k,\infty)$. Two
elements of $GL(k,\infty)$ lie in the same orbit if and
only if they have the same image in $\Bbb R^{\infty}$. 
We have a surjective mapping 
$\pi:GL(k,\infty)\to G(k,\infty)$,   
given by $\pi(A)= A(\Bbb R^k)$, where the inverse images of points
are exactly the $GL(k)$-orbits.

Similarly we have
$$O(k,\infty)=\{A\in L(\Bbb R^k,\Bbb R^{\infty}): A^t\o A= Id_k\}.$$
Now the Lie group $O(k)$ of all isometries of $\Bbb R^k$ acts
freely from the right on the space $O(k,\infty)$. Two
elements of $O(k,\infty)$ lie in the same orbit if and
only if they have the same image in $\Bbb R^{\infty}$. 
The projection $\pi:GL(k,\infty)\to G(k,\infty)$ 
restricts to a surjective mapping
$\pi:O(k,\infty)\to G(k,\infty)$
and the inverse images of points
are now exactly the $O(k)$-orbits.

\proclaim{\nmb.{6.4}. Lemma (Iwasawa decomposition)} 
Let $T(k;\Bbb R)$ be the group of all upper triangular $k\x
k$-matrices with positive 
entries on the main diagonal. Then each $B\in GL(k,\infty)$ can
be written in the form
$B= p(B)\o q(B)$, with unique
$p(B)\in O(k,\infty)$ and $q(B)\in T(k)$. The mapping
$q:GL(k,\infty)\to T(k)$ is real analytic.
\endproclaim

\proclaim{\nmb.{6.5}. Theorem} The following are a real
analytic principal fiber bundles:
$$\gather (O(k,\infty;\Bbb R),\pi,G(k,\infty;\Bbb R),O(k,\Bbb R)),\\
(GL(k,\infty;\Bbb R),\pi,G(k,\infty;\Bbb R),GL(k,\Bbb R)),\\
(GL(k,\infty;\Bbb R),p,O(k,\infty;\Bbb R),T(k;\Bbb R)).
\endgather$$
The last one is trivial. 

The embeddings $\Bbb R^n\to \Bbb R^\infty$ induce real analytic
embeddings, which respect the principal right actions of all the
structure groups 
$$\gather O(k,n)\to O(k,\infty),\\
GL(k,n)\to GL(k,\infty),\\
G(k,n)\to G(k,\infty).
\endgather$$
All these cones are inductive limits in the category of real
analytic (and smooth) manifolds.
\endproclaim

\proclaim{\nmb.{6.6}. Theorem} The following manifolds are
real analytically diffeomorphic to the homogeneous spaces indicated:
$$\gather 
GL(k,\infty)\cong {GL(\infty)}\left/
    {\pmatrix Id_k & L(\Bbb R^k,\Bbb R^{\infty-k})\\
        0 & GL(\infty-k)\endpmatrix}\right. \\
O(k,\infty)\cong O(\infty)/Id_k\x O(\infty-k)\\
G(k,\infty)\cong O(\infty)/O(k)\x O(\infty-k).
\endgather$$
The {universal vector bundle}
$(E(k,\infty),\pi,G(k,\infty),\Bbb R^k)$ is defined as the
associated bundle
$$\align E(k,\infty)&=O(k,\infty)[\Bbb R^k] \\
&=\{(Q,x): x\in Q\}\subset G(k,\infty)\x \Bbb R^\infty.
\endalign$$
The tangent bundle of the Grassmannian is then given by 
$$TG(k,\infty) = L(E(k,\infty),E(k,\infty)^\bot).$$
\endproclaim

\proclaim{\nmb.{6.7} Theorem} The
principal bundle $(O(k,\infty),\pi,G(k,\infty))$ is classifying
for finite dimensional principal $O(k)$-bundles and carries a 
universal real analytic $O(k)$-connection
$\om\in\Om^1(O(k,\infty),\frak o(k))$.

This means: For each finite dimensional smooth or real analytic principal
$O(k)$-bundle $P\to M$ with principal connection $\om_P$ there
is a smooth or real analytic mapping $f:M\to G(k,\infty)$ such
that the pullback $O(k)$-bundle $f^*O(k,\infty)$ is isomorphic to $P$ and the
pullback connection $f^*\om = \om_P$ via this isomorphism.
\endproclaim

\subheading{\nmb.{6.8}. The Lie group $GL(\infty,\Bbb R)$} The
canonical embeddings $\Bbb R^n\to \Bbb R^{n+1}$ onto the first
$n$ coordinates induce injections $GL(n)\to GL(n+1)$. The
inductive limit is 
$$GL(\infty):= \varinjlim_nGL(n)$$
in the category of sets. Since each $GL(n)$ also injects into
$L(\Bbb R^\infty,\Bbb R^\infty)$ we can visualize $GL(\infty)$
as the set of all $\Bbb N\x \Bbb N$-matrices which are
invertible and differ from the identity in finitely many entries
only. 

We also consider the Lie algebra $\frak g\frak l(\infty)$
of all $\Bbb N\x \Bbb N$-matrices with only finitely many
nonzero entries, which is isomorphic to 
$\Bbb R^{(\Bbb N\x\Bbb N)}$, 
and we equip it with this convenient vector space structure.
Then $\frak g\frak l(\infty)
=\varinjlim_n\frak g\frak l(n)$ in the category of
real analytic mappings.

\proclaim{Claim} $\frak g\frak l(\infty)=L(\Bbb R^{\Bbb
N},\Bbb R^{(\Bbb N)})$ as convenient vector spaces. Composition
is a bounded bilinear mapping on $\frak g\frak l(\infty)$.
\endproclaim

\proclaim{\nmb.{6.9}. Theorem} $GL(\infty)$ is a real analytic
Lie group modeled on $\Bbb R^\infty$, with Lie algebra 
$\frak g\frak l(\infty)$ and is the inductive limit of the Lie
groups $GL(n)$ in the category of real analytic manifolds. 
The exponential mapping 
is well defined, is real analytic and a local real analytic diffeomorphism
onto a neighborhood of the identity. The
Campbell-Baker-Hausdorff formula gives a real analytic mapping
near 0 and expresses the multiplication on $GL(\infty)$ via $\exp$.
The determinant $\det:GL(\infty)\to \Bbb R\setminus 0$ is a real analytic
homomorphism. We have a real analytic left action
$GL(\infty)\x\Bbb R^\infty\to\Bbb R^\infty$, such that $\Bbb
R^\infty\setminus 0$ is one orbit, but the injection 
$GL(\infty)\hookrightarrow L(\Bbb R^\infty,\Bbb R^\infty)$
does not generate the topology.
\endproclaim
\demo{Proof} Since the exponential mappings are compatible with
the inductive limits all these assertions follow from the
inductive limit property.
\qed\enddemo

\proclaim{\nmb.{6.10}. Theorem} Let $\frak g$ be a Lie
subalgebra of $\frak g\frak l(\infty)$. Then there is a smoothly
arcwise connected
splitting Lie subgroup $G$ of $GL(\infty)$ whose Lie algebra is
$\frak g$. The exponential mapping of $GL(\infty)$ restricts to
that of $G$, which is local diffeomorphism near zero. The
Campbell Baker Hausdorff formula gives a real analytic mapping
near 0 and has the usual properties, also on $G$.
\endproclaim
\demo{Proof} Let $\frak g_n:= \frak g\cap \frak g\frak l(n)$, a
finite dimensional Lie subalgebra of $\frak g$. 
Then $\bigcup \frak g_n = \frak g$. The convenient structure 
$\frak g =\varinjlim_{n}\frak g_n$ coincides with the structure
inherited as a complemented subspace, since $\frak g\frak l(\infty)$ carries
the finest locally convex structure.

So for each $n$ there is a connected Lie subgroup $G_n\subset
GL(n)$ with Lie algebra $\frak g_n$. Since 
$\frak g_n\subset \frak g_{n+1}$ we have $G_n\subset G_{n+1}$
and we may consider $G:= \bigcup_n G_n\subset GL(\infty)$. 
Each $g\in G$ lies in some $G_n$ and may be connected to $Id$
via a smooth curve there, which is also smooth curve in $G$, so
$G$ is smoothly arcwise connected. 

All mappings $\exp|\frak g_n:\frak g_n\to G_n$ are local real
analytic diffeomorphisms near 0, so $\exp:\frak g\to G$ is also
a local real analytic diffeomorphism near zero onto an open
neighborhood of the identity in $G$. The rest is clear.  
\qed\enddemo

\subheading\nofrills{\bf\nmb.{6.11}. Examples.\newline 
The Lie group $SO(\infty,\Bbb R)$\usualspace} is the
inductive limit 
$$SO(\infty):= \varinjlim_nSO(n)\subset GL(\infty).$$
It is the connected Lie subgroup of $GL(\infty)$ with the Lie algebra
$\frak o(\infty) = \{X\in \frak g\frak l(\infty):X^t=-X\}$ of
skew elements. 
Obviously we have 
$$\multline 
SO(\infty)=\{A\in GL(\infty):\langle Ax,Ay\rangle =\langle x,y\rangle \\
\text{ for all }x,y\in \Bbb R^\infty\text{ and }\det(A)=1\}.
\endmultline$$

\subheading\nofrills{\bf The Lie group\usualspace }  $O(\infty)$ is the
inductive limit
$$\multline O(\infty):= \varinjlim_nO(n)\subset GL(\infty)\\
    =\{A\in GL(\infty):\langle Ax,Ay\rangle=\langle x,y\rangle
    \text{ for all }x,y\in \Bbb R^\infty\}.
\endmultline$$
It has two connected components, that of the identity is
$SO(\infty)$. 

\subheading\nofrills{\bf The Lie group\usualspace } $SL(\infty)$ is the
inductive limit
$$\align SL(\infty):&= \varinjlim_nSL(n)\subset GL(\infty)\\
    &= \{A\in GL(\infty): \det(A)=1\}.
\endalign$$
It is the connected Lie subgroup with Lie algebra
$\frak s\frak l(\infty)=\{X\in \frak g\frak l(\infty):
\operatorname{Trace}(X)=0\}$.

\subheading{\nmb.{6.12}} We stop here to give examples. Of course 
this method is also applicable for the complex versions of the most 
important homogeneous spaces. This will be treated elsewhere.

\enddocument